\documentclass{article}

\usepackage[utf8]{inputenc}
\usepackage{titling}
\usepackage{lipsum}
\usepackage[space]{grffile}
\usepackage{graphicx}
\usepackage{url}
\usepackage{amsmath}
\usepackage{amssymb}
\usepackage{hyperref}
\usepackage{float}

\allowdisplaybreaks

\usepackage{verbatim}

\usepackage{xcolor}

\usepackage{pbox}

\newcommand{\ginv}[1]{{#1}^{\mbox{\tiny -U}}}
\newcommand{\pinv}[1]{{#1}^{\mbox{\tiny -P}}}

\newcommand{\inv}[1]{{#1}^{\mbox{\tiny -1}}}
\newcommand{\sinv}[1]{{#1}^{\mbox{\tiny -S}}}

\begin{document}

\begin{Large}
\begin{center}
{\bf Triple-Block Generalized Inverses for Control}\\
{\bf Applications with Mixed Consistency Requirements}
\end{center}
\end{Large}
\vspace{-11pt}
\begin{center} 
{\large Jeffrey Uhlmann}\\
Dept.\ of Electrical Engineering and Computer Science\\
University of Missouri - Columbia
\end{center}
\vspace{0pt}
\begin{quote}
\begin{center}
{\bf Abstract}
\end{center}
\vspace{-3pt}
Previous work on block-partitioned mixed generalized 
inverses is extended from two subsets of system variables
with distinct consistency requirements to three
subsets. Does not include any significant 
theoretical contributions.
\end{quote}

\vspace{3pt}

{\bf INTRODUCTION}\vspace{3pt}

Previous work has identified the need to select appropriate 
generalized inverses to enforce critical application-specific 
assumptions, such as whether system behavior should be 
invariant with respect to rigid rotations or with respect to 
changes of units on state variables \cite{uhlmann}. For example, 
the familiar Moore-Penrose (MP) pseudoinverse is only applicable
when invariance is needed with respect to orthogonal/unitary 
transformations of the input and/or output of the system. The 
unit-consistent (UC) generalized inverse, by contrast, is only applicable
when invariance is needed with respect to diagonal transformations
of the input and/or output of the system, e.g., as determined by
choices of units on state variables. 

It is not widely recognized that the choice of generalized inverse is
critical for robust performance of a control system. Consequently,
it is common to see control-system implementations that simply
replace all matrix inverses with the MP inverse (e.g., replace {\em inv}
with {\em pinv} when using Matlab) to defend against singular
states, such as arise during gimball lock or other loss-of-rank
configurations. The fact is that there are an infinite number of
generalized inverses, each of which sacrifices a specific subset of
properties of the true matrix inverse \cite{ben,uhlmann}. 

The term {\em consistency} is used to describe the behavior of a
system with respect to a set of transformations \cite{jkusim,uhlmann0}. 
For example, if the behavior of a robotic system is expected to be invariant
up to arbitrary orthogonal transformations of an assumed Euclidean
coordinate frame, then the MP inverse is appropriate. In other words,
the behavior of the system should be identical up to a rotation of the 
coordinate frame. If, however, the MP inverse is applied with
variables involving arbitrary choices of units, e.g., centimeters
versus meters for lengths, then the system will exhibit different
behavior depending on that choices of units, which is clearly
undesirable \cite{boz1,boz2,ucrga,rafal}. 

A challenge that arises in most nontrivial systems is that there
exist distinct subsets of variables with different consistency
requirements. For example, one subset may be defined in 
arbitrary units of length, whereas another may be defined in
an arbitrarily-oriented Euclidean coordinate frame. This
means that reliable control (i.e., insensitive to the choice of
units or rotation of coordinate frame) must involve a generalized
inverse that treats the different subsets of variables in 
different ways. This can be expressed in terms of a block
partition of the system matrix $M$ with $m$ variables requiring
unit consistency and the remaining $n$ variables requiring 
rotation consistency:
\begin{eqnarray*}
   \begin{array}{rcl}
      M & = &
         \left[ \begin{array}{cc} { W} & { X}\\
                             { Y}& { Z}\end{array} \right] \vspace{-4pt}
                       \begin{array}{l} {\}}~m \\{\}}~n \end{array}\\
                        & &  
                                \begin{array}{c}    \hspace{-2pt}                                                     
                                   \;\underbrace{\;}_m\;\underbrace{\;}_n
                                \end{array}
   \end{array} . 
\end{eqnarray*}
It has been shown \cite{uhlmann0,uhlmann} that the mixed inverse can be obtained from this block-partitioned form as
\begin{equation}
   M^{\mbox{\tiny -J}} = 
   \left[
     \begin{array}{cc}
     \ginv{({ W}-{ X}\pinv{{ Z}}{ Y})} & -\ginv{{ W}}{ X}\pinv{({ Z}-{ Y}\ginv{{ W}}{ X})} \\
      -\pinv{{ Z}}{ Y}\ginv{({ W}-{ X}\pinv{{ Z}}{ Y})} & \pinv{({ Z}-{ Y}\ginv{{ W}}{ X})} 
     \end{array}
    \right]  \label{dualBlock}
\end{equation}
where superscript -U denotes the UC inverse, the superscript -P denotes the
MP pseudoinverse, and -J denotes the specialized joint inverse of the system
matrix. This dual-block {\em mixed inverse} form has been applied in multiple robotic
control systems \cite{boz1,boz2}, and it is sufficient for most applications. 
In the following section we generalize the mixed-inverse structure to also include
a third distinct generalized inverse. In other words, we provide a solution
for the case in which there are {\em three} subsets of variables, each of which
has distinct consistency requirements.

~\\
{\bf THE TRIPLE-BLOCK MIXED INVERSE}\vspace{3pt}

The triple-block mixed inverse demands a formulation of the generalized
inverse of a matrix of the following form:
\begin{equation}
M ~=~ \left[\begin{matrix}
\;R & S & T\, \\
\;U & V & W\, \\
\;X & Y & Z\,
\end{matrix}\right] \nonumber
\end{equation}

The dual-block solution of Eq.\,(1) can be concisely expressed. 
Unfortunately, the complexity of the solution for the $3\times 3$ block matrix
is considerably greater\footnote{It does not appear that any of the current
computer algebra systems, e.g., Mathematica and Maple, are capable of
performing the kinds of symbolic block-matrix manipulations needed to
automate the solution of problems of this kind. These expressions were 
therefore derived laboriously by hand, but have not been empirically 
verified. To avoid introducing transcription mistakes, no factoring 
of expressions has been performed.}. The following are expressions 
for each block of M$^{\mbox{\tiny -J}}$ in unsimplified form with superscripts
\{-a,-b,-c\} denoting three distinct generalized inverses:
\begin{eqnarray*}
\mbox{\bf M}^{\mbox{\tiny -J}}_{11} & = &  
((R^{-a} S ((-U R^{-a} S+V)^{-b}+(-U R^{-a} S+V)^{-b} (-U R^{-a} T+W) \\
&~& (-(-X R^{-a} S+Y) (-U R^{-a} S+V)^{-b} (-U R^{-a} T+W)-X R^{-a} T+Z)^{-c}\\
&~& (-X R^{-a} S+Y) (-U R^{-a} S+V)^{-b})-R^{-a} T (-(-X R^{-a} S+Y)\\
&~& (-U R^{-a} S+V)^{-b} (-U R^{-a} T+W)-X R^{-a} T+Z)^{-c} (-X R^{-a} S+Y)\\
&~& (-U R^{-a} S+V)^{-b}) U+(-R^{-a} S (-U R^{-a} S+V)^{-b} (-U R^{-a} T+W)\\
&~& (-(-X R^{-a} S+Y) (-U R^{-a} S+V)^{-b} (-U R^{-a} T+W)-X R^{-a} T+Z)^{-c}+\\
&~& R^{-a} T (-(-X R^{-a} S+Y) (-U R^{-a} S+V)^{-b} (-U R^{-a} T+W)\\
&~& -X R^{-a} T+Z)^{-c}) X) R^{-a}+R^{-a}\\
&~&\\ 
\mbox{\bf M}^{\mbox{\tiny -J}}_{12}  &= &
-R^{-a} S ((-U R^{-a} S+V)^{-b}+(-U R^{-a} S+V)^{-b} (-U R^{-a} T+W)\\
&~& (-(-X R^{-a} S+Y) (-U R^{-a} S+V)^{-b}(-U R^{-a} T+W)-X R^{-a} T+Z)^{-c}\\
&~& (-X R^{-a} S+Y) (-U R^{-a} S+V)^{-b})+R^{-a} T (-(-X R^{-a} S+Y) (-U R^{-a}S+V)^{-b}\\
&~& (-U R^{-a} T+W)-X R^{-a} T+Z)^{-c} (-X R^{-a} S+Y) (-U R^{-a} S+V)^{-b}   \\
&~&\\
\mbox{\bf M}^{\mbox{\tiny -J}}_{13}   & = &
R^{-a} S (-U R^{-a} S+V)^{-b} (-U R^{-a} T+W) (-(-X R^{-a} S+Y)\\
&~& (-U R^{-a} S+V)^{-b} (-U R^{-a} T+W)-X R^{-a} T+Z)^{-c}-R^{-a} T\\
&~& (-(-X R^{-a} S+Y) (-U R^{-a} S+V)^{-b} (-U R^{-a} T+W)-X R^{-a} T+Z)^{-c} \\
&~&\\
\mbox{\bf M}^{\mbox{\tiny -J}}_{21}   & = &
-(((-U R^{-a} S+V)^{-b}+(-U R^{-a} S+V)^{-b} (-U R^{-a} T+W)\\
&~& (-(-X R^{-a} S+Y) (-U R^{-a} S+V)^{-b} (-U R^{-a} T+W)-X R^{-a} T+Z)^{-c}\\
&~& (-X R^{-a} S+Y) (-U R^{-a} S+V)^{-b}) U-(-U R^{-a} S+V)^{-b} (-U R^{-a} T+W)\\
&~& (-(-X R^{-a} S+Y) (-U R^{-a} S+V)^{-b} (-U R^{-a} T+W)-X R^{-a} T+Z)^{-c} X) R^{-a} \\
&~&\\
\mbox{\bf M}^{\mbox{\tiny -J}}_{22}   & = &
(-U R^{-a} S+V)^{-b}+(-U R^{-a} S+V)^{-b} (-U R^{-a} T+W)\\
&~& (-(-X R^{-a} S+Y) (-U R^{-a} S+V)^{-b} (-U R^{-a} T+W)-X R^{-a} T+Z)^{-c}\\
&~& (-X R^{-a} S+Y) (-U R^{-a} S+V)^{-b} \\
&~&\\
\mbox{\bf M}^{\mbox{\tiny -J}}_{23}  & = &
-(-U R^{-a} S+V)^{-b} (-U R^{-a} T+W) (-(-X R^{-a} S+Y)\\
&~& (-U R^{-a} S+V)^{-b} (-U R^{-a} T+W)-X R^{-a} T+Z)^{-c} \\
&~&\\
\mbox{\bf M}^{\mbox{\tiny -J}}_{31}   & = &
-(-(-(-X R^{-a} S+Y) (-U R^{-a} S+V)^{-b} (-U R^{-a} T+W)\\
&~& -X R^{-a} T+Z)^{-c} (-X R^{-a} S+Y) (-U R^{-a} S+V)^{-b} U+(-(-X R^{-a} S+Y)\\
&~& (-U R^{-a} S+V)^{-b} (-U R^{-a} T+W)-X R^{-a} T+Z)^{-c} X) R^{-a} \\
&~&\\
\mbox{\bf M}^{\mbox{\tiny -J}}_{32}  & = &
-(-(-X R^{-a} S+Y) (-U R^{-a} S+V)^{-b} (-U R^{-a} T+W)\\
&~& -X R^{-a} T+Z)^{-c} (-X R^{-a} S+Y) (-U R^{-a} S+V)^{-b} \\
&~&\\
\mbox{\bf M}^{\mbox{\tiny -J}}_{33}  & = &
(-(-X R^{-a} S+Y) (-U R^{-a} S+V)^{-b}\\
&~& (-U R^{-a} T+W)-X R^{-a} T+Z)^{-c} 
\end{eqnarray*}

For purposes of redundancy, and to facilitate practical
implementation, the following is a conversion from the
same derivation as above into functional form. The
three user-selected generalized inverses are denoted as
\texttt{InvA}, \texttt{InvB}, and \texttt{InvC}.
Indentation and square brackets for functions are
purely for ease of parsing.

\begin{quote}
\begin{verbatim}
Mi11 = 
(
  (InvA[R]*S*
    (
      InvB[-U*InvA[R]*S+V]+Inverse[-U*InvA[R]*S+V]*
      (-U*InvA[R]*T+W)*
      Inverse[-
        (-X*InvA[R]*S+Y)*InvB[-U*InvA[R]*S+V]*
        (-U*InvA[R]*T+W)-X*InvA[R]*T+Z
      ]*(-X*InvA[R]*S+Y)*InvB[-U*InvA[R]*S+V]
    )-InvA[R]*T*
    InvC[-
      (-X*InvA[R]*S+Y)*
      InvB[-U*InvA[R]*S+V]*(-U*InvA[R]*T+W)- 
      X*InvA[R]*T+Z
    ]*(-X*InvA[R]*S+Y)*InvB[-U*InvA[R]*S+V]
  )*U+
  (-
    InvA[R]*S*InvB[-U*InvA[R]*S+V]*
    (-U*InvA[R]*T+W)*
    Inverse[-
      (-X*InvA[R]*S+Y)*InvB[-U*InvA[R]*S+V]*
      (-U*InvA[R]*T+W)-X*InvA[R]*T+Z
    ]+InvA[R]*T*
    InvC[-
      (-X*InvA[R]*S+Y)*InvB[-U*InvA[R]*S+V]*
      (-U*InvA[R]*T+W)-X*InvA[R]*T+Z
    ]
  )*X
)*InvA[R]+InvA[R] ;

Mi12 = -
InvA[R]*S*
(
  InvB[-U*InvA[R]*S+V]+InvB[-U*InvA[R]*S+V]*
  (-U*InvA[R]*T+W)*
  InvC[-
    (-X*InvA[R]*S+Y)*InvB[-U*InvA[R]*S+V]*
    (-U*InvA[R]*T+W)-X*InvA[R]*T+Z
  ]*(-X*InvA[R]*S+Y)*
  InvB[-U*InvA[R]*S+V])+InvA[R]*T*
  InvC[-
    (-X*InvA[R]*S+Y)*InvB[-U*InvA[R]*S+V]*
    (-U*InvA[R]*T+W)-X*InvA[R]*T+Z
  ]*(-X*InvA[R]*S+Y)*InvB[-U*InvA[R]*S+V] ;

Mi13 = 
InvA[R]*S*InvB[-U*InvA[R]*S+V]*
(-U*InvA[R]*T+W)*
Inverse[-
  (-X*InvA[R]*S+Y)*InvB[-U*InvA[R]*S+V]*
  (-U*InvA[R]*T+W)-X*InvA[R]*T+Z
]-InvA[R]*T*
Inverse[-
  (-X*InvA[R]*S+Y)*InvB[-U*InvA[R]*S+V]*
  (-U*InvA[R]*T+W)-X*InvA[R]*T+Z
] ;

Mi21 = -
(
  (
    Inverse[-U*InvA[R]*S+V]+InvB[-U*InvA[R]*S+V]*
    (-U*InvA[R]*T+W)*
    InvC[-
      (-X*InvA[R]*S+Y)*InvB[-U*InvA[R]*S+V]*
      (-U*InvA[R]*T+W)-X*InvA[R]*T+Z
    ]*(-X*InvA[R]*S+Y)*InvB[-U*InvA[R]*S+V]
  )*U-InvB[-U*InvA[R]*S+V]*(-U*InvA[R]*T+W)*
  InvC[-
    (-X*InvA[R]*S+Y)*InvB[-U*InvA[R]*S+V]*
    (-U*InvA[R]*T+W)-X*InvA[R]*T+Z
  ]*X
)*InvA[R] ;

Mi22 = 
InvB[-U*InvA[R]*S+V]+InvB[-U*InvA[R]*S+V]*
(-U*InvA[R]*T+W)*
InvC[-
  (-X*InvA[R]*S+Y)*InvB[-U*InvA[R]*S+V]*
  (-U*InvA[R]*T+W)-X*InvA[R]*T+Z
]*(-X*InvA[R]*S+Y)*InvB[-U*InvA[R]*S+V] ;

Mi23 = -
InvB[-U*InvA[R]*S+V]*(-U*InvA[R]*T+W)*
InvC[-
  (-X*InvA[R]*S+Y)*Inverse[-U*InvA[R]*S+V]*
  (-U*InvA[R]*T+W)-X*InvA[R]*T+Z
] ;

Mi31 = - 
(-
  InvC[-
    (-X*InvA[R]*S+Y)*InvB[-U*InvA[R]*S+V]*
    (-U*InvA[R]*T+W)-X*InvA[R]*T+Z
  ]*(-X*InvA[R]*S+Y)*InvC[-U*InvA[R]*S+V]*U+
  Inverse[-
    (-X*InvA[R]*S+Y)*InvB[-U*InvA[R]*S+V]*
    (-U*InvA[R]*T+W)-X*InvA[R]*T+Z
  ]*X
)*InvA[R] ;

Mi32 = - 
Inverse[-
  (-X*InvA[R]*S+Y)*InvB[-U*InvA[R]*S+V]*
  (-U*InvA[R]*T+W)-X*InvA[R]*T+Z
]*(-X*InvA[R]*S+Y)*InvB[-U*InvA[R]*S+V] ;

Mi33 = 
Inverse[-
  (-X*InvA[R]*S+Y)*InvB[-U*InvA[R]*S+V]*
  (-U*InvA[R]*T+W)-X*InvA[R]*T+Z
] ;
\end{verbatim}
\end{quote}

The formatting of the above reveals several
small and large subexpressions that can be factored 
out for efficiency. However, this explicit form can
serve as a baseline for verifying optimized
variants.

~\\
{\bf ILLUSTRATIVE EXAMPLE}\vspace{3pt}

The MP inverse provides consistency with respect to left
or right unitary transforms of the system matrix. The 
UC inverse provides consistency with respect to left or
right nonsingular diagonal transforms of the system matrx.
These two are by far the most common forms of 
required consistency that arise in practical applications.

Another class of tranforms are more structured with
respect to the kind of consistency demanded of the 
left and right transformations. One is 
similarity, i.e., requiring consistency
with respect to transformations of the system
matrix $M$ as $SM\inv{S}$ for arbitrary 
nonsingular matrix $S$. 

A similarity-consistent (SC) generalized inverse
can be obtained by decomposing $M$ as
\begin{equation}
M ~=~ FC\inv{F}, \nonumber
\end{equation}
where $C$ is the Frobenius Canonical Form (FCF)
of $M$, which then permits the SC inverse to
be evaluated using the MP inverse of $C$ as
\begin{equation}
\sinv{M} ~=~ F\pinv{C}\inv{F}. \nonumber
\end{equation}
This form is theoretically tractable to compute,
unlike the Jordan Normal Form (JNF), to provide
similarity consistency\footnote{The Drazin inverse
can also be used in cases in which it produces
a result with the same rank as $M$ \cite{uhlmann}.
There are also various speculative approaches \cite{jkusim}. 
One class involves the scaling 
of $M$ by a positive value large enough relative to 
the precision of the entries that they can be 
treated as integers to permit the Frobenius
Rational Canonical Form or Smith Normal Form
to be used, both of which can be efficiently
computed. However, the numerical implications
of integerising $M$ are not well understood.}\cite{jkusim}.

Now, given a system with $m$ variables requiring
Euclidean rotational consistency; $n$ variables 
requiring diagonal unit consistency; and $p$
variables requiring similarity consistency:
\begin{eqnarray*}
   \begin{array}{rcl}
      M & = &
         \!\left[ \begin{array}{ccc} { R} & { \,\;S} & { \;T}\\
                             { U}& { \,\;V} &{ \;W}\\
                             { X}& { \,\;Y} &{ \;Z}\end{array} \right] \vspace{-4pt}
                       \begin{array}{l} {\}}~m \\{\}}~n \\{\}}~p\end{array}\\
                        & &  
                                \begin{array}{c}    \hspace{-2pt}                                                     
                                   \;\underbrace{\;}_m\;\underbrace{\;}_n\;\underbrace{\;}_p
                                \end{array}
   \end{array} 
\end{eqnarray*}
the appropriate generalized matrix inverse of $M$ can
be constructed using the formulation of the previous 
section with \texttt{InvA} provided in the form of the
MP inverse; \texttt{InvB} provided in the form of the UC
inverse; and \texttt{InvC} provided in the form of the
SC inverse.

~\\
{\bf DISCUSSION}\vspace{3pt}

We have provided an explicit
formulation of the triple-block mixed
generalized matrix inverse. This complements
the dual-block formulation for cases in which
three subsets of system variables have
distinct consistency requirements. 

The triple-block case is in fact redundant because the general
solution can be recursively constructed from repeated 
applications of the dual-block solution. Specifically, given 
$k$ sets of variables $S_1...S_k$, with corresponding sizes 
$n_1...n_k$, a mixed inverse for sets $S_1$ and $S_2$
can be constructed using the dual-block solution. 
With this, their respective state variables can be treated 
as a single set of size $n_1\!+\!n_2$ with consistency 
conditions satisfied by the newly-constructed mixed 
inverse. In other words, the number of sets of variables 
can be reduced by 1. This can then be repeated to the 
base case of two sets.
 
The motivation for deriving an explicit 
triple-block solution is purely for convenience.  
Specifically, while the dual-block inverse is 
much more likely to find use in practical 
control applications, we anticipate a likelihood 
that increasing awareness of consistency 
issues will lead to recognition of cases 
in which the triple-block formulation is 
needed. Having an explicit plug-and-play
solution will hopefully reduce the perceived
effort to apply that recognition.

Future work will examine consistency 
considerations arising in physics models, e.g., 
dynamical systems for which it is necessary to 
rigorously establish that key consistency/invariance 
properties are preserved across reduced-rank 
manifolds (or where a covariant derivative 
becomes singular).

\bibliographystyle{plain}

\begin{thebibliography}{11}


\bibitem{rafal}
Rafal Qasim Al Yousuf and Jeffrey Uhlmann, ``On Use of the Moore-Penrose Pseudoinverse for Evaluating the RGA of Non-Square Systems,'' {\em Iraqi Journal of Computers, Communications, Control \& Systems Engineering} (IJCCCE), 3:21, pp.\ 89-97, DOI: 10.33103/uot.ijccce.21.3.8, 2021.

\bibitem{ben}
Ben-Israel, Adi; Greville, Thomas N.E., ``Generalized inverses: Theory and applications, " 
{\em New York: Springer}, 2003.

\bibitem{uhlmann}
Jeffrey Uhlmann, ``A Generalized Matrix Inverse that is Consistent with Respect to Diagonal Transformations,''
{\em SIAM Journal on Matrix Analysis} (SIMAX), Vol.\ 39:3, 2018.

\bibitem{uhlmann0}
Jeffrey Uhlmann,  ``Unit Consistency, Generalized Inverses, and
Effective System Design Methods,'' arXiv:1604.08476v2 [cs.NA] 11 Jul 2017 (2015).

\bibitem{jkusim}
Jeffrey Uhlmann, ``A Rank-Preserving Generalized Matrix Inverse for
Consistency with Respect to Similarity,'' arXiv:1804.07334v1 [cs.NA] 19 Apr 2018.

\bibitem{ucrga}
Jeffrey Uhlmann, ``On the Relative Gain Array (RGA) with Singular and
Rectangular Matrices,'' {\em Applied Mathematics Letters}, Vol.\ 93, 2019.

\bibitem{boz1}
Bo Zhang and Jeffrey Uhlmann, ``Applying a Unit-Consistent Generalized Matrix Inverse for Stable Control of Robotic Systems,'' {\em ASME J. of Mechanisms and  Robotics}, 11(3), 2019.

\bibitem{boz2}
Bo Zhang and Jeffrey Uhlmann, ``Examining a Mixed Inverse Approach for Stable Control of a Rover,'' {\em International Journal of Control Systems and Robotics}, 5, 1-7, 2020.


\end{thebibliography}

\end{document}